\documentclass[12pt]{article}
\usepackage{amssymb,amsfonts,amsmath, psfrag}
\newtheorem{theo}{Theorem}

\newtheorem{pro}[theo]{Proposition}

\def\qed{\hfill \rule{4pt}{7pt}}
\def\MA{\mathbb{A}}
\begin{document}
\begin{center}
{\bf \Large $q$-Identities from Lagrange and Newton Interpolation}
\end{center}
\vspace{0.3cm}
\begin{center}
  {\bf Amy Mei Fu}\\
  Center for Combinatorics \\
  Nankai University, Tianjin 300071, P.R. China\\
  Email: fmu@eyou.com
 \end{center}
\begin{center}
 {\bf Alain Lascoux}\\
 Nankai University, Tianjin 300071, P.R. China\\
 Email: Alain.Lascoux@univ-mlv.fr\\
 CNRS, IGM Universit\'e  de
 Marne-la-Vall\'ee\\
 77454 Marne-la-Vall\'ee Cedex, France\\
 \end{center}

 \noindent {\bf ABSTRACT:}  Combining Newton and Lagrange
interpolation, we give $q$-identities which  generalize results of
Van Hamme, Uchimura, Dilcher and Prodinger.
\section{introduction}

Van Hamme \cite{Van} gave the following identity involving Gauss
polynomials, see also Andrews \cite{And},
  \begin{equation}\label{L1}
    \sum_{i=1}^{n}{n \brack i}\frac{(-1)^{i-1}q^{\binom{i+1}{2}}}{1-q^i}=
    \sum_{i=1}^n\frac{q^i}{1-q^i},
  \end{equation}
 where ${n \brack i}$ is the Gauss polynomials defined by${n \brack i}=(q;q)_{n}\left((q;q)_{i}(q;q)_{n-i}\right)^{-1}$
 with $(z;q)_{n}=(1-z)(1-zq)\cdots (1-zq^{n-1})$.

Uchimura \cite{Uch} generalize (\ref{L1}) as following:
\begin{equation}\label{L2}
\sum_{i=1}^n{n \brack
i}\frac{(-1)^{i-1}q^{\binom{i+1}{2}}}{1-q^{i+m}}=\sum_{i=1}^n\frac{q^i}{1-q^i}\left/{i+m
\brack i}, \ m\geq 0\right.
\end{equation}
and by Dilcher \cite{Dil}:
\begin{equation}\label{L3}
\sum_{1 \leq i \leq n}{n \brack
i}(-1)^{i-1}\frac{q^{\binom{i}{2}+mi}}{(1-q^i)^m}=\sum_{1\leq
i_{1}\leq i_{2}\leq \cdots \leq i_{m}\leq
n}\frac{q^{i_{1}}}{1-q^{i_{1}}}\cdots
\frac{q^{i_{m}}}{1-q^{i_{m}}}.
\end{equation}

Prodinger \cite{Pro}  mentioned the following identity as a
$q$-analogue of Kirchenhofer's \cite{Kirl} formula,
\begin{eqnarray}\label{L17}
\sum_{i=0, i\neq M}{n \brack
i}\frac{(-1)^{i-1}q^{\binom{i+1}{2}}}{1-q^{i-M}}=(-1)^Mq^{\binom{M+1}{2}}{n
\brack M }\sum_{i=0, i\neq M}\frac{q^{i-M}}{1-q^{i-M}}
\end{eqnarray}
and  explained how to obtain all these formulas by using Cauchy
residues.

We shall show that in fact, all the above formulas are a direct
consequence of Newton  and Lagrange interpolation.

Given two finite sets of variables  $\MA$
 and $\mathbb{B}$, we denote by $R(\MA,\mathbb{B})$ the product
 $\prod_{a \in \MA, b \in \mathbb{B}}(a-b)$, and by  $\MA
\setminus \mathbb{B} $  the set difference of $\MA$ and
$\mathbb{B}$.

Let $\MA=\{x_{1},x_{2},\cdots  \}$, $\MA_{n}=\{x_{1},x_{2},\cdots,
x_{n}\}$, for any $n \geq 0$.  Lagrange wrote the following
summation:
\begin{equation}\label{L4}
f(x) = \sum_{i=1}^{n}f(x_{i})\frac{R(x,\MA\setminus
x_{i})}{R(x_{i}, \MA\setminus x_{i})},\  {\rm mod} R(x, \MA_{n}).
\end{equation}

On the other hand, Newton's development is:
\begin{equation}\label{Add5}
f(x)=f(x_{1})+f\partial_{1} R(x,
\MA_{1})+f\partial_{1}\partial_{2}R(x, \MA_{2})+\cdots,
\end{equation}
where $\partial_{i}$, $i \geq 1$, operating on its left, is
defined by
\begin{eqnarray*}
f(x_{1},x_{2},\cdots, x_{i},x_{i+1},\cdots)\partial_{i}
=\frac{f(\cdots, x_{i},x_{i+1},\cdots)-f(\cdots,
x_{i+1},x_{i},\cdots)}{x_{i}-x_{i+1}}.
\end{eqnarray*}
Taking $x_{n+1}=x$, Newton's summation terminates:
\begin{eqnarray*}
f(x)=f(x_{1})+f\partial_{1} R(x,
\MA_{1})+f\partial_{1}\partial_{2}R(x,
\MA_{2})+\cdots+f\partial_{1}\cdots\partial_{n}R(x,\MA_{n}).
\end{eqnarray*}

Newton's and Lagrange's expressions imply the same remainder
$f\partial_{1}\cdots\partial_{n}R(x,\MA_{n})$, and  the polynomial
$g_{n}(x)=f(x)-f\partial_{1}\cdots\partial_{n}R(x,\MA_{n})$ is the
only polynomial of degree $\leq n-1$ such that $f(x_{i})=g(x_{i}),
\ 1 \leq i \leq n$.

Taking now  $f(x)=(y-x)^{-1}$, since
\begin{eqnarray*}
\frac{1}{y-x_{1}}\partial_{1} \cdots
\partial_{n-1}=\frac{1}{(y-x_{1})\cdots(y-x_{n})},
\end{eqnarray*}
 one gets, by comparing Newton's and Lagrange's expressions, the identity:
\begin{equation}\label{L14}
\sum_{i=0}^{n-1}\frac{R(x,\MA_{i})}{R(y,\MA_{i+1})}
=\sum_{i=1}^n\frac{f(x_{i})R(x,\MA \setminus x_{i})}{R(x_{i},
\MA\setminus x_{i}
)}=\frac{1}{y-x}-\frac{R(x,\MA_{n})}{R(y,\MA_{n})(y-x)}.
\end{equation}

Letting $x=1$, we derive:
\begin{equation}\label{L16}
\sum_{i=0}^{n-1}\frac{R(1,\MA_{i})}{R(y,\MA_{i+1})}
=\sum_{i=1}^n\frac{f(x_{i})R(1,\MA\setminus x_{i})}{R(x_{i},
\MA\setminus x_{i} )}.
\end{equation}

Expanding $(y-x)^{-1}=\sum_{i}y^{-i-1}x^{i}$ and taking the
coefficient of $y^{-m-1}$, we have:
\begin{equation}\label{L5}
x_{1}^{m}\partial_{1}\cdots\partial_{n-1}=\sum_{i=1}^n
\frac{x_{i}^m}{\prod_{j \neq
i}(x_{i}-x_{j})}=h_{m-n+1}(x_{1},x_{2},\cdots ,x_{n}).
\end{equation}
Recall that complete functions \cite{Mac} $h_{k}$ are defined by
\begin{eqnarray*}
h_{k}(x_{1},x_{2},\cdots ,x_{n})=\sum_{1\leq i_{1}\leq i_{2}\leq
\cdots \leq i_{k}\leq n}x_{i_{1}}x_{i_{2}}\cdots x_{i_{k}}.
\end{eqnarray*}

\section{ General identities}

The identities that we present just corresponds to taking
$\mathbb{A}=\{\frac{a-bq}{c-zq},\frac{a-bq^2}{c-zq^2},\cdots
 \}$ in Lagrange or Newton interpolation. In that case, the
 products $R(x_i , \MA \setminus x_i)$ are immediate to write.

\begin{pro}\label{Pro1}
 Let $m, n \in \mathbb{N}$, $\tau=m-n+1$, $a, b, z, q $ be variables ,
 and  $\MA=\{\frac{a-bq}{c-zq},\frac{a-bq^2}{c-zq^2},
 \cdots,\frac{a-bq^n}{c-zq^n} \}$.
 We have:
 \begin{eqnarray}\label{L8}
 \lefteqn{\sum_{1\leq i_{1}\leq i_{2}\leq \cdots \leq i_{\tau}\leq
 n}\frac{a-bq^{i_{1}}}{c-zq^{i_{1}}}\frac{a-bq^{i_{2}}}{c-zq^{i_{2}}}\cdots
 \frac{a-bq^{i_{\tau}}}{c-zq^{i_{\tau}}}=} \nonumber \\
  &&\frac{c^n(zq/c;q)_{n}}{(q;q)_{n}(az-bc)^{n-1}}\sum_{i=1}^n{n \brack
 i}\frac{(-1)^{i-1}
 q^{\binom{i+1}{2}-ni}(1-q^i)(a-bq^i)^m}{(c-zq^i)^{\tau+1}}. \nonumber
 \\
\end{eqnarray}
 In particular, for $m=n$:
\begin{eqnarray}\label{L9}
\lefteqn{ \sum_{i=1}^n\frac{a-bq^i}{c-zq^i}=}\nonumber \\
&&\frac{c^n(zq/c;q)_{n}}{(q;q)_{n}(az-bc)^{n-1}}\sum_{i=1}^n{n
\brack
 i}\frac{(-1)^{i-1}
q^{\binom{i+1}{2}-ni}(1-q^i)(a-bq^i)^n}{(c-zq^i)^2}\nonumber\\
\ \
\end{eqnarray}

\end{pro}

{\bf Proof.} For our choice of $\MA$, then $R(x_i, \MA\setminus
x_i )=\prod_{j \neq
 i}\frac{(az-bc)(q^i-q^j)}{(c-zq^i)(c-zq^j)}$, and
\begin{eqnarray*}
\lefteqn{\sum_{i=1}^n\frac{((a-bq^i)/(c-zq^i))^m}{\prod_{j \neq
 i}(\frac{a-bq^i}{c-zq^i}-\frac{a-bq^j}{c-zq^j})}}\\
 &=&\frac{c^n(zq/c;q)_{n}}{(q;q)_{n}(az-bc)^{n-1}}\sum_{i=1}^n {n
\brack i} \frac{(-1)^{i-1}
 q^{\binom{i+1}{2}-ni}(1-q^i)(a-bq^i)^m}{(c-zq^i)^{m-n+2}}.
\end{eqnarray*}
\qed

Consider the case $a=0$, $b=-1$, $c=1$ and $z=1$, i.e.
$\MA=\{\frac{q}{1-q},\frac{q^2}{1-q^2}, \cdots
\frac{q^n}{1-q^n}\}$.

Eq (\ref{L9}) implies Van Hamme identity :
\begin{eqnarray*}
\sum_{i=1}^n\frac{q^i}{1-q^i}
&=&\frac{(q;q)_{n}}{(q;q)_{n}}\sum_{i=1}^n{n \brack
i}\frac{(-1)^{i-1}q^{\binom{i+1}{2}-ni}(1-q^i)q^{ni}}{(1-q^i)^2}\\
&=&\sum_{i=1}^n {n \brack
i}\frac{(-1)^{i-1}q^{\binom{i+1}{2}}}{1-q^i}.
\end{eqnarray*}

From (\ref{L8}), we have,
\begin{eqnarray*}
 \sum_{1\leq i_{1}\leq i_{2}\leq \cdots \leq i_{m-n+1}\leq
n}\frac{q^{i_{1}}}{1-q^{i_{1}}}\frac{q^{i_{2}}}{1-q^{i_{2}}}\cdots
\frac{q^{i_{m-n+1}}}{1-q^{i_{m-n+1}}}=\sum_{i=1}^n{n \brack
i}\frac{(-1)^{i+1}q^{\binom{i+1}{2}-ni}q^{mi}}{(1-q^i)^{m-n+1}},
\end{eqnarray*}
replacing $m-n+1$ by $m$,  we derive (\ref{L3}).

Take now $a=1$, $b=0$, $c=0$ , $z=-1$, then
$\MA=\{q^{-1},q^{-2},\cdots, q^{-n}\}$, and Eq.(\ref{L16})
implies:
\begin{equation*}
R.H.S=\sum_{i=0}^{n-1}\frac{R(1, \MA_{i})}{R(y,
\MA_{i+1})}=\sum_{i=0}^{n-1}\frac{\prod_{j=1}^i(1-q^{-j})}{\prod_{j=1}^{i+1}(y-q^{-j})}
=-\sum_{i=0}^{n-1}\frac{q^{i+1}(q;q)_{i}}{(yq;q)_{i+1}},
\end{equation*}
\begin{multline*}
L.H.S=\sum_{i=1}^n\frac{f(x_{i})R(1, \MA\setminus
x_{i})}{R(x_{i},\MA\setminus x_{i}
)}\\
=\sum_{i=1}^n\frac{(y-q^{-i})^{-1}\prod_{j \neq
i}(1-q^{-j})}{\prod_{j \neq i}(q^{-i}-q^{-j})}
 =\sum_{i=1}^n{n
\brack i}\frac{(-1)^i q^{\binom{i+1}{2}}}{1-yq^i}.
\end{multline*}
Clearly:
\begin{eqnarray*}
\sum_{i=1}^n{n \brack i}\frac{(-1)^{i-1}
q^{\binom{i+1}{2}}}{1-yq^i}
=\sum_{i=0}^{n-1}\frac{q^{i+1}(q;q)_{i}}{(yq;q)_{i+1}}=\sum_{i=1}^{n}\frac{q^{i}(q;q)_{i-1}}{(yq;q)_{i}}.
\end{eqnarray*}
In the special case $y=q^m$, we obtain (\ref{L2}).

 We  take now $\MA=\{q^{M},q^{M-1}, q^{M-2}, \cdots , q,
q^{-1}, \cdots q^{M-n}\}$, with $M \in \mathbb{N}$, and $y=1$.
 Eq(\ref{L16}) becomes:
\begin{equation*}
 R.H.S=\sum_{i=0}^{n-1}\frac{R(1,\MA_{i})}{R(1,\MA_{i+1})}=-\sum_{i=0
, i \neq M}^n\frac{q^{i-M}}{1-q^{i-M}},
\end{equation*}
and
\begin{multline*}
L.H.S = \sum_{i=0,i\neq M}^n\frac{1}{1-q^{M-i}}\frac{\prod_{j\neq
i,M}(1-q^{M-j})}{\prod_{j \neq i,M}(q^{M-i}-q^{M-j})}\\
=\sum_{i=0 i\neq
M}^n\frac{(-1)^{M+i}q^{\binom{i+1}{2}-\binom{M+1}{2}}(q;q)_{M}(q;q)_{n-M}}{(q;q)_{i}(q;q)_{n-i}},
\end{multline*}
which proves (\ref{L17}).

 As a final comment, we would like to stress that we have
 just used simple alphabets in Newton and Lagrange interpolation; it is easy to generalize the above formulas by taking
 more sophisticated alphabets.

\end{document}